\documentclass[a4paper]{jpconf} \usepackage{graphicx}
\usepackage{color} \usepackage[latin1]{inputenc}
\usepackage[english]{babel}
\usepackage{amsfonts,amssymb,amsmath,latexsym,upref}
\usepackage{graphics,amssymb,amsthm,amsmath}
\newtheorem{Thm}{Theorem}[section]
\newtheorem{Def}[Thm]{Definition} 
\newtheorem{Lem}[Thm]{Lemma} 

\newtheorem{Prop}[Thm]{Proposition} 

\newcommand{\qdet}{{\operatorname{det}}_q}
\newcommand{\qDet}{{\operatorname{Det}}_q}
\newenvironment{dedication} { \itshape 
\raggedleft } {\par 
   
} \begin{document} \title{Special classes of homomorphisms
between generalized Verma modules for ${\mathcal U}_q(su(n,n))$}

\author{Hans Plesner Jakobsen}

\address{Department of Mathematical Sciences, University of
Copenhagen, Denmark}

\ead{jakobsen@math.ku.dk}

\begin{abstract}We study homomorphisms between quantized
generalized Verma modules
$M(V_{\Lambda})\stackrel{\phi_{\Lambda,\Lambda_1}}{\rightarrow}M(V_{\Lambda_1})$
for ${\mathcal U}_q(su(n,n))$. There is a natural notion of
degree for such maps, and if the map is of degree $k$, we write
$\phi^k_{\Lambda,\Lambda_1}$. We examine when one can have a
series of such homomorphisms
$\phi^1_{\Lambda_{n-1},\Lambda_{n}}\circ
\phi^1_{\Lambda_{n-2},\Lambda_{n-1}}\circ\cdots\circ\phi^1_{\Lambda,\Lambda_1}
=\textrm{Det}_q$, where $\textrm{Det}_q$ denotes the map $M(V_{\Lambda})\ni p\rightarrow \textrm{Det}_q\cdot p\in M(V_{\Lambda_n})$. If, classically, $su(n,n)^{\mathbb
C}={\mathfrak p}^-\oplus(su(n)\oplus su(n)\oplus {\mathbb
C})\oplus {\mathfrak p}^+$, then
$\Lambda=(\Lambda_L,\Lambda_R,\lambda)$ and
$\Lambda_n=(\Lambda_L,\Lambda_R,\lambda+2)$. The answer is then that
$\Lambda$ must be one-sided in the sense that either
$\Lambda_L=0$ or $\Lambda_R=0$ (non-exclusively). There are
further demands on $\lambda$ if we insist on ${\mathcal
U}_q({\mathfrak g}^{\mathbb C})$ homomorphisms. However, it is
also interesting to loosen this to considering only ${\mathcal
U}^-_q({\mathfrak g}^{\mathbb C})$ homomorphisms, in which case
the conditions on $\lambda$ disappear. \smallskip By duality,
there result have implications on covariant quantized
differential operators. We finish by giving an explicit, though
sketched, determination of the full set of ${\mathcal
U}_q({\mathfrak g}^{\mathbb C})$ homomorphisms
$\phi^1_{\Lambda,\Lambda_1}$.

\end{abstract}

\begin{dedication} Dedicated to I.E. Segal (1918-1998) in
commemoration of the centenary of his birth. \end{dedication}

\section{Introduction} Generalized and quantized Verma modules
have physically attractive properties similar to the Fock space.
There is a ``vacuum vector'', here called a highest weight
vector, which is annihilated by the ``upper diagonal''
operators, is an eigenvector for the ``diagonal operators'', and
which generate the whole space when acted upon by the algebra of
``lower diagonal operators''. Since it may happen that there is
a second vacuum vector, it is of interest to determine cases in
which this may happen. This is further interesting because by
duality, such cases correspond to quantized covariant
differential operators such as the Maxwell equations. We give
here a complete proof of the one-sidedness and we give a sketch
of the case of an arbitrary first order. Further details as well
as the dual picture will appear in a forthcoming article. For
the ``classical'' analogue, see e.g. \cite{jak}. On a personal
note: The explicitness presented here is in line with how
mathematical physics was taught to me by Segal, my Ph.D.
advisor. \bigskip

\section{Set-up}
\begin{equation}\begin{array}{ccccccccc}{\mathfrak g}^{\mathbb
C}&=&su(n,n)^{\mathbb C}\qquad&=&{\mathfrak k}^{\mathbb C}\oplus
{\mathfrak p}&=&{\mathfrak p}^-\oplus {\mathfrak k}^{\mathbb
C}\oplus {\mathfrak p}^+&=&{\mathfrak p}^-\oplus {\mathfrak
k}^{\mathbb C}\oplus {\mathfrak p}^+, \\{\mathfrak k}^{\mathbb
C}&=&su(n)^{\mathbb C}\oplus {\mathbb C}\oplus su(n)^{\mathbb
C}&=&{\mathfrak k}_L^{\mathbb C}\oplus {\zeta} \oplus {\mathfrak
k}_R^{\mathbb C}\ .&&&& \end{array} \end{equation}
\begin{equation}\begin{array}{ccccccccc} \zeta \textrm{ is the
center, }&{\mathfrak p}^{\pm}\textrm{ are abelian ${\mathcal
U}({\mathfrak k}^{\mathbb C})$ modules,}&\textrm{ and }&
{\mathcal U}({\mathfrak g}^{\mathbb C})&={\mathcal P}({\mathfrak
p}^-)\cdot {\mathcal U}({\mathfrak k}^{\mathbb C})\cdot
{\mathcal P}({\mathfrak p}^+) \end{array} .\end{equation}

We let the simple roots be denoted
$\Pi=\{\mu_1,\dots,\mu_{n-1}\}\cup\{\beta\}\cup
\{\mu_1,\dots,\mu_{n-1}\}$, where $\beta$ is the unique
non-compact roots and where the decomposition of simple roots
corresponds to the decomposition of ${\mathfrak k}^{\mathbb C}$
above. \smallskip In the quantum group ${\mathcal
U}_q({\mathfrak g}^{\mathbb C})$, we denote the generators by
$E_\alpha,F_\alpha,K_\alpha^{\pm1}$ for $\alpha\in\Pi$. There
are also decompositions \begin{eqnarray} {\mathcal
U}_q({\mathfrak g}^{\mathbb C})&=&{\mathcal A}_q^-\cdot
{\mathcal U}_q({\mathfrak k}^{\mathbb C})\cdot {\mathcal
A_q}^+,\\ {\mathcal U}_q({\mathfrak k}^{\mathbb C})&=&{\mathcal
U}_q({\mathfrak k}_L^{\mathbb C})\cdot {\mathbb
C}[K_\beta^{\pm1}] \cdot {\mathcal U}_q({\mathfrak k}_R^{\mathbb
C}). \end{eqnarray} Here, ${\mathcal A}_q^\pm$ are quadratic
algebras which are furthermore ${\mathcal U}_q({\mathfrak
k}^{\mathbb C})$. Specifically, \begin{eqnarray} {\mathcal
A}_q^-&=&{\mathbb C}[W_{i,j}\mid i,j=1,\dots,n\},\\ {\mathcal
A}_q^+&=&{\mathbb C}[Z_{i,j}\mid i,j=1,\dots,n\}, \end{eqnarray}
with relations \begin{eqnarray}\label{a}Z_{ij}Z_{ik} &=&
q^{-1}Z_{ik}Z_{ij} \textrm{ if }j < k;\\ \label{b}Z_{ij}Z_{kj}
&=& q^{-1}Z_{kj}Z_{ij}\textrm{ if }i< k;\\\label{c} Z_{ij}Z_{st}
&=& Z_{st}Z_{ij} \textrm{ if }i < s\textrm{ and }t <
j;\\\label{cross} Z_{ij}Z_{st} &=&
Z_{st}Z_{ij}-(q-q^{-1})Z_{it}Z_{sj} = \textrm{ if }i < s\textrm{
and } j < t. \end{eqnarray} The algebra ${\mathcal A}_q^-$ have
the same relations, but the algebras ${\mathcal A}_q^\pm$ are
different as ${\mathcal U}_q({\mathfrak k}^{\mathbb C})$
modules. The elements $Z_{ij}$ and $W_{ij}$ are constructed by
means of the Lusztig operators. References \cite{l} and
\cite{jan} are general references of much of this. Using the
Serre relations one gets, setting $\mu_0=Id$, \begin{Lem}
\label{2.1}\begin{eqnarray}Z_{i,j}&=&T_{\nu_{j-1}}
T_{\nu_{j-2}}\dots T_{\nu_{0}}\cdot T_{\mu_{i-1}}
T_{\mu_{i-2}}\dots T_{\mu_{0}}(E_\beta),\\
W_{i,j}&=&T_{\nu_{j-1}} T_{\nu_{j-2}}\dots T_{\nu_{0}}\cdot
T_{\mu_{i-1}} T_{\mu_{i-2}}\dots
T_{\mu_{0}}(F_\beta).\label{26}\end{eqnarray} \end{Lem}

\smallskip

For later use, we give the relations in the full algebra:

\begin{eqnarray} E_{\mu_k} W_{i,j}&=&W_{i,j}E_{\mu_k}\textrm{ if
}k\neq i-1,\\
E_{\mu_k}W_{i,j}^a&=&(-q)[a]W_{i-1,j}W_{i,j}^{a-1}K_{\mu_k}+W^a_{i,j}E_{\mu_k}
\textrm{ if }k=i-1,\\ F_{\mu_k}
W_{i,j}&=&W_{i,j}F_{\mu_k}\textrm{ if }k\neq i,i-1,\\
F_{\mu_k}W_{i,j}^a&=&-q^{-1}[a]W_{i,j}^{a-1}W_{i+1,j}+q^{-a}W^a_{i,j}F_{\mu_k}
\textrm{ if }k=i,\\ F_{\mu_k}
W_{i,j}&=&qW_{i,j}F_{\mu_k}\textrm{ if }k= i-1,\\ F_{\mu_k}
Z_{i,j}&=&Z_{i,j}F_{\mu_k}\textrm{ if }k\neq i-1,\\
F_{\mu_k}Z_{i,j}^a&=&[a]Z_{i-1,j}Z_{i,j}^{a-1}K^{-1}_{\mu_k}+W^a_{i,j}E_{\mu_k}
\textrm{ if }k=i-1,\\ E_{\mu_k}
Z_{i,j}&=&Z_{i,j}E_{\mu_k}\textrm{ if }k\neq i,i-1,\\
E_{\mu_k}Z_{i,j}^a&=&[a]Z_{i,j}^{a-1}Z_{i+1,j}+q^{-a}Z^a_{i,j}E_{\mu_k}
\textrm{ if }k=i,\\ E_{\mu_k}
Z_{i,j}&=&qZ_{i,j}E_{\mu_k}\textrm{ if }k= i-1. \end{eqnarray}

There are similar formulas for the commutators involving
$E_{\nu_k}$ and $F_{\nu_k}$. If e.g. $S$ denotes the obvious
automorphism defined on generators by $W_{ij}\rightarrow
W_{j,i}$, and similarly, $Z_{ij}\rightarrow Z_{j,i}$ then
$E_{\nu_k}=SE_{\mu_k}S$ and $F_{\nu_k}=SF_{\mu_k}S$.

\section{Finite dimensional ${\mathcal U}_q({\mathfrak
k}^{\mathbb C})$ modules} A non-zero vector $v_\Lambda$ of a
finite dimensional module $V_\Lambda$ of ${\mathcal
U}_q({\mathfrak k}^{\mathbb C})$ is a highest weight vector of
highest weight $\Lambda$, and $V_\Lambda$ is a highest weight
module of highest weight $\lambda$,
if\begin{equation}\begin{array}{lccccccc} \forall i=1,\dots,n-1:
K_{\mu_i}^{\pm1}=q^{\pm\lambda^\mu_i}v_\Lambda,&
K_{\nu_i}^{\pm1}= q^{\pm \lambda^\nu_i}v_\Lambda, &\textrm{ and
} K_{\beta}^{\pm1}=q^{\pm\lambda}v_\Lambda.\\ \textrm{Finally,
}{\mathcal U}^+_q({\mathfrak k}^{\mathbb C})v_\Lambda=0,
&\textrm{and } {\mathcal U}^-_q({\mathfrak k}^{\mathbb
C})v_\Lambda=V . \end{array}\end{equation} We set
$\Lambda=((\lambda^\mu_1,\dots,\lambda^\mu_{n-1}),(\lambda^\nu_1,\dots,
\lambda^\nu_{n-1});\lambda)=(\Lambda_L, \Lambda_R,\lambda)$.

As a vector space, $V_\Lambda=V_{\Lambda_L}\otimes
V_{\Lambda_R}$ where $V_{\Lambda_L}$ and $V_{\Lambda_R}$ are
highest weight representations of ${\mathcal U}_q({\mathfrak
k}_L^{\mathbb C})$ and ${\mathcal U}_q({\mathfrak k}_R^{\mathbb
C})$, respectively, of highest weights
$\Lambda_L=(\lambda^\mu_1,\dots,\lambda^\mu_{n-1})$ and
$\Lambda_R=(\lambda^\nu_1,\dots,\lambda^\nu_{n-1})$,
respectively. The highest weight vector can then be written as
$v_\Lambda=v_{\Lambda_L}\otimes v_{\Lambda_R}$ with the
stipulation that $K_\beta^{\pm1}v_{\Lambda_L}\otimes
v_{\Lambda_R}=q^{\pm\lambda}v_{\Lambda_L} \otimes
v_{\Lambda_R}$.

\bigskip

\section{Generalized quantized Verma modules and their
homomorphisms}

Consider a finite dimensional module
$V_\Lambda=V_{\Lambda_L,\Lambda_R,\lambda}$ over ${\mathcal
U}_q({\mathfrak k}^{\mathbb C})$ with highest weight is defined
by $\Lambda=(\Lambda_L,\Lambda_R,\lambda)$ where
$\Lambda_L=(\lambda^\mu_1,\lambda^\mu_2,\dots,\lambda^\mu_{n-1},0)$,
$\Lambda_R=(\lambda^\nu_1,\lambda^\nu_2,\dots,\lambda^\nu_{n-1},0)$,
and $\lambda\in{\mathbb C}$.

We extend such a module to a ${\mathcal U}_q({\mathfrak
k}^{\mathbb C}){\mathcal A}_q^+$ module, by the same name, by
letting ${\mathcal A}_q^+$ act trivially.

\begin{Def} The quantized generalized Verma module
$M(V_\Lambda)$ is given by
\begin{equation}M(V_\Lambda)={\mathcal U}_q({\mathfrak
g}^{\mathbb C})\bigotimes_{{\mathcal U}_q({\mathfrak k}^{\mathbb
C}){\mathcal A}_q^+}V_\Lambda \end{equation} with the natural
action from the left. \end{Def} \smallskip

As a vector space, \begin{equation}M(V_\Lambda)={\mathcal
A}_q^-\otimes V_\Lambda. \end{equation} We are interested in
structure preserving homomorphisms between quantized generalized
Verma modules. We call such maps intertwiners, covariants, or
equivariants, indiscriminately. Dually, they will be quantized
covariant differential operators. In abstract notation, the
structure under investigation is\begin{equation}Hom_{{\mathcal
U}_q({\mathfrak g}^{\mathbb
C})}(M(V_{\Lambda}),M(V_{\Lambda_1})). \end{equation}However,
for the time being we will
consider\begin{equation}\label{former}Hom_{{\mathcal
A}_q^-{\mathcal U}_q({\mathfrak k}^{\mathbb
C})}(M(V_{\Lambda}),M(V_{\Lambda_1})). \end{equation} An element
$\phi_{\Lambda,\Lambda_1}$ in the latter space is completely
determined by the ${\mathcal U}_q({\mathfrak k}^{\mathbb C})$
equivariant map, denoted by the same symbol:
\begin{equation}V_{\Lambda}\stackrel
{\phi_{\Lambda,\Lambda_1}}{\rightarrow} {\mathcal A}_q^-\otimes
V_{\Lambda_1}\textrm{ leads to } {\mathcal A}_q^-\otimes
V_{\Lambda}\stackrel{\phi_{\Lambda,\Lambda_1}}{\rightarrow}
{\mathcal A}_q^-\otimes V_{\Lambda_1}.
\end{equation}Specifically, ${\phi_{\Lambda,\Lambda_1}}$ does
not depend on $\lambda$ and is completely given by the condition
that the image of the highest weight vector
${\phi_{\Lambda,\Lambda_1}}(v_{\Lambda})$ is a highest weight
vector for ${\mathcal U}_q({\mathfrak k}^{\mathbb C})$. For the
map $\phi_{\Lambda,\Lambda_1}$ to belong to the former space
(\ref{former}) it is necessary, and sufficient that,
additionally, ($Z_\beta$ acting in $M(V_{\Lambda_1}$))
\begin{equation}\label{zbeta} Z_\beta\left(
{\phi_{\Lambda,\Lambda_1}}(v_{\Lambda})\right)=0. \end{equation}
This equation depends heavily on $\lambda$. It is clear that
such maps, whether of the first or second kind, can be combined:
\begin{equation} \phi_{\Lambda_1,\Lambda_2}\circ
\phi_{\Lambda,\Lambda_1}=\phi_{\Lambda,\Lambda_2} \end{equation}
though it may happen that the composite is zero.

We use the terminology of degree of elements of ${\mathcal
A}_q^-$ in the obvious way, and we let, for $k=1,\dots$, $
{\mathcal A}_q^-(k)$ denote the ${\mathcal U}_q({\mathfrak
k}^{\mathbb C})$ module spanned by homogeneous elements of
degree $k$. If the elements $p_{ij}$ all belong to ${\mathcal
A}_q^-(k)$, we write $\phi^k_{\Lambda,\Lambda_1}$.

\medskip

{\bf General Problem:} When is it possible to have
$\phi^1_{\Lambda_{n-1},\Lambda_{n}}\circ
\phi^1_{\Lambda_{n-2},\Lambda_{n-1}}\circ\cdots\circ\phi^1_{\Lambda,\Lambda_1}
=\qDet$? In this case, if
$\Lambda=(\lambda_L,\Lambda_R,\lambda)$ , then
$\Lambda_{n}=(\Lambda_L,\Lambda_R,\lambda+2)$.

\bigskip

\section{Laplace expansion}

If $m=n$, one may define the quantum determinant $det_q$ in
${\mathcal A}_q-$ as follows: \begin{eqnarray}\label{2b}
\qdet(n)=\qdet&=&\Sigma_{\sigma\in
S_n}(-q^{-1})^{\ell(\sigma)}W_{1,\sigma(1)}W_{2,\sigma(2)}
\cdots W_{n,\sigma(n)}\\\label{3b}&=&\Sigma_{\delta\in
S_n}(-q^{-1})^{\ell(\delta)}W_{\delta(1),1}W_{\delta(2),2}
\cdots W_{\delta(n),n}. \end{eqnarray}

\medskip

If $m=n$ and $I=\{i_1<1_2< \dots <
i_{n-1}\}=\{1,2,\cdots,n\}\setminus\{i\},J=\{j_1<j_2< \dots <
j_{n-1}\}=\{1,2,\cdots,n\}\setminus\{j\}$, we set
\begin{eqnarray}\label{2} A(i,j)&=&\Sigma_{\sigma\in
S_{n-1}}(-q^{-1})^{\ell(\sigma)}W_{i_1,j_{\sigma(1)}}W_{i_2,j_{\sigma(2)}}
\cdots
W_{i_{n-1},j_{\sigma({n-1})}}\\\label{3}&=&\Sigma_{\tau\in
S_{n-1}}(-q^{-1})^{\ell(\tau)}
W_{i_{\tau(1)},j_1}W_{i_{\tau(2)},j_2} \cdots
W_{i_{\tau({n-1})},j_{n-1}}.\end{eqnarray} These elements are
quantum $(n-1)\times(n-1)$ minors. The following was proved by
Parshall and Wang \cite{pw}: \begin{Prop}$\qdet$ is central.
Furthermore, let $i,k\le n$ be fixed integers. Then

\begin{eqnarray}\label{325}\delta_{i,k}\qdet&=&\sum_{j=1}^n(-q^{-1})^{j-k}W_{i,j}
A(k, j)=\sum_j(-q^{-1})^{i-j} A(i,j)W_{k,j}\\\label{326}
&=&\sum_j(-q^{-1})^{j-k}W_{j,i}A(j,k) =
\sum_j(-q^{-1})^{i-j}A(j,i)W_{j,k}.\end{eqnarray} \end{Prop}

\bigskip

\section{1. order}

Any finite dimensional highest weight representation of
${\mathcal U}_q({\mathfrak k}^\mathbb C)$ of the form
$\Lambda=(\Lambda_L,\Lambda_R,\lambda)$ in which either
$\Lambda_L=0$ or $\Lambda_R=0$ will be called one-sided. We will
now give an explicit form for a highest weight vector $v_1$ of
an irreducible sub-representation of ${\mathcal A}^-(1)\otimes
V_{\Lambda=(\Lambda_L,0,\lambda)}$. Specifically,
consider\begin{equation}\label{hwv}
v_1=W_{N+1,1}v_0+W_{N,1}u_Nv_0
+W_{N-1,1}u_{N-1,1}v_0+W_{N-2,1}u_{N-2}v_0
+\dots+W_{1,1}u_{1}v_0, \end{equation} where $\forall
i=1,\dots,N: u_i\in {\mathcal
U}^{-\mu_N-\mu_{N-1}-\dots-\mu_{i}}_q({\mathfrak k}^\mathbb C)$.
\smallskip

Because of this, we first want to consider a basis of ${\mathcal
U}_q^{-\mu_N+\cdots- \mu_\ell}({\mathfrak k}^{\mathbb C}_L)$.

Set ${\mathcal E}_{\ell,N}=\{\ell,\ell+1,\cdots,N\}\subseteq
\{1,2,\dots,n-1\}$. Any sequence
$I_{\ell,N}=(i_\ell,i_{\ell+1},\cdots,i_N)$ made up of pairwise
different elements of ${\mathcal E}^\mu_{\ell,N}$ defines a
non-zero element \begin{equation}
F_{\mu_{i_\ell}}F_{\mu_{i_{\ell+1}}}\cdots
F_{\mu_{i_N}}=F^\mu(I_{\ell,N})\in {\mathcal
U}_q^{-\mu_\ell+\cdots- \mu_N}({\mathfrak k}^{\mathbb C}_L).
\end{equation} We will call such a sequence {\bf allowed}. We
reserve the name $E_{\ell,N}$ for the special sequence
$(\ell,\ell+1,\cdots,N)$. \medskip

We will say that a transposition $(i_\ell,i_{\ell+1},\dots,
i_k,i_{k+1},\dots,i_N)\rightarrow (i_\ell,i_{\ell+1},\dots,
i_{k+1},i_k,\dots,i_N)$ is legal if $\mid i_{k+1}-i_k\mid>1$.

Recall that $F_{\mu_i}F_{\mu_j}=F_{\mu_j}F_{\mu_i}$ if $\vert
i-j\vert>1$. We will say that two allowed sequences $I^{(1)}$
and $I^{(2)}$ are equivalent if one can be obtained from the
other by a series of legal transpositions. It is clear that any
allowed sequence $I$ can be brought, uniquely, and by legal
transpositions, into the form $J_1J_2\cdots J_r$ which is the
concatenation of sequences $J_t$ that are either descending or
ascending, and such that the following are satisfied: Firstly,
the elements of $J_s$ are smaller than the elements of $J_t$ if
$s<t$, and $\cup_s J_s=\{\ell,\ell+1,\dots,N\}$. Secondly, two
neighboring sequences cannot both be ascending (maximality), and
thirdly, singletons are ascending.

We denote by ${\mathcal J}_{\ell,N}$ the set of such sequences.
\smallskip The following is then obvious:
\begin{Prop}\begin{equation} \{F^\mu(I_{\ell,N})\mid
I_{\ell,N}\in {\mathcal J}_{\ell,N}\} \end{equation} is a basis
of $ {\mathcal U}_q^{-\mu_\ell+\cdots- \mu_N}({\mathfrak
k}^{\mathbb C}_L)$. \end{Prop}

We furthermore have from e.g. \cite[lemma~6.27]{dpw}:

\begin{Prop}\label{basis} Let $V=V(\Lambda_L)$ be a finite
dimensional highest weight representation of ${\mathcal
U}_q({\mathfrak k}^{\mathbb C}_L)$ with
$\Lambda_L=(\lambda^\mu_1,\lambda^\mu_2,\cdots,
\lambda^\mu_{n-1})$ satisfying: $\lambda^\mu_\ell>0,
\lambda^\mu_{\ell+1}>0,\dots, \lambda^\mu_N>0$. Let $v_0$ denote
a highest weight vector (unique up to a non zero constant). Then
\begin{equation}\{F^\mu(I_{\ell,N})v_0\mid I_{\ell,N}\in
{\mathcal J}_{\ell,N}\} \end{equation} is a basis of
$V^{\Lambda_L-\mu_\ell+\cdots- \mu_N}$. \end{Prop}

If $I_{\ell,N}=J_1J_2\cdots J_s\in {\mathcal J}_{\ell,N}$ as
above, we attach to it a sequence
$C^\mu(I_{\ell,N})=(c_{i_\ell},c_{i_{ell+1}},\cdots c_{i_N})$
where $c_k=a_k$ if either $i_k$ belongs to an ascending
sub-sequence $J_x$ of $I_{\ell,N}$ or if $i_k$ is the biggest
element in a descending sub-sequence $J_y$ of $I_{\ell,N}$.
Here, $x,y\in\{1, 2,\dots,s\}$. In the remaining cases,
$c_{i_k}=b_{i_k}$. We furthermore set
$f^\mu(C^\mu(I_{\ell,N}))=\prod_{t=\ell}^Nc_{i_t}$.

\smallskip

We can then state, maintaining the assumptions from
Lemma~\ref{basis}:

\begin{Prop}\label{Prop-basis} If the vector $v_1$ in
(\ref{hwv}) is a highest weight vector in ${\mathcal
A}_1^-\otimes V(\Lambda_L)$ then \begin{equation}\forall
\ell=1,\dots,N:u_\ell=\sum_{I_{\ell,N}\in{\mathcal
J}_{\ell,N}}f^\mu(C^\mu(I_{\ell,N}))F^\mu(I_{\ell,N})v_0.
\end{equation} \end{Prop} Later, we shall find it convenient to
set ${\mathcal J}_{N+1,N}=\emptyset$ and
$f^\mu(C^\mu(\emptyset))=1 =F^\mu(\emptyset)$. Likewise,
$E_{N+1,N}=\emptyset$.

\medskip

Our general case of interest is where we only assume
$\lambda^\mu_N\neq0$. Bear in mind that in the sequence
$C(I_{\ell,N})$, $c_0=b_i$ signals that the corresponding
$\mu_i$, taking part in $F(I_{\ell,N})$, can be moved all the
way to the right without changing $F(I_{\ell,N})$. If we allow
$\lambda^\mu_i=0$ this means that such elements, when applied to
$v_0$, give zero. Hence if we let ${\mathcal
Z}_{\ell,N}=\{i={\ell,\cdots,N}\mid\lambda_i^\mu=0\}$ and if we
let ${\mathcal J}_{\ell,N}^{\mathcal Z}$ denote those sequences
$I$ in ${\mathcal J}_{\ell,N}$ for which any index $i$ from
${\mathcal Z}_{\ell,N}$ either belongs to an increasing sequence
or is the biggest index in a decreasing sequence, then we have:
\begin{Prop}
\label{6.4}\begin{equation}\{F^\mu(I_{\ell,N})v_0\mid
I_{\ell,N}\in {\mathcal J}^{\mathcal Z}_{\ell,N}\}
\end{equation} is a basis of $V^{\Lambda_L-\mu_\ell+\cdots-
\mu_N}$. \end{Prop}

\smallskip

Clearly there is an analogue to Proposition~\ref{Prop-basis} for
this general case (just as long as $\lambda^\mu_N>0$).

\smallskip

There is yet another helpful way to view the various sets
${\mathcal J}_{\ell,N}$, $\ell=N,N-1,\dots,1$, namely as a
labeled, directed rooted tree with root at $F_{\mu_N}$:
\begin{equation} \begin{array}{rcl} &F^\mu(I_{\ell,N})&\\
\stackrel{L_{\ell-1}}{\swarrow}&&\stackrel{R_{\ell-1}}{\searrow}\\
F_{\mu_{\ell-1}}F^\mu(I_{\ell,N})
&&F^\mu(I_{\ell,N})F_{\mu_{\ell-1}} \end{array}. \end{equation}
Here, it is really only the relative positions of $F_{\mu_\ell}$
and $F_{\mu_{\ell-1}}$ that matter.

If we have $\lambda_i^\mu=0$ we just modify the tree by removing
all branches labeled by $R_i$ - as well as everything above
these branches - from the tree. (In this picture, the root is
lowest.)

In this way, there is an obvious bijection between the paths in
the modified tree and the basis.

\medskip

We now return to (\ref{hwv}). To obtain the following equations,
it is used that $E_{\mu_{i-1}}(W_{i,j})=-qW_{i-1,j}K_{\mu_{i-1}}
+(W_{i,j})E_{\mu_{i-1}}$, which follows from Lemma~\ref{2.1}.
Furthermore, for the vector in (\ref{hwv}) to be a ${\mathcal
U}_q({\mathfrak k}^\mathbb C)$ highest weight vector we clearly
only need to look at ${\mathcal U}_q({\mathfrak k}_L^\mathbb
C)$. Here we must have: \begin{eqnarray}\forall i=1,\dots,N:
(-q)W_{i,1}K_{\mu_i}u_{i+1}v_0
+W_{i,1}E_{\mu_i}u_{i}v_0&=&0\\\forall i,j=1,\dots,N:
E_{\mu_j}u_{i}v_0&=&0\textrm{ if }i\neq j. \end{eqnarray} We
assume throughout that $\lambda^\mu_N\neq0$.

Using Proposition~\ref{Prop-basis}, we set $u_{N+1}=1$ and
\begin{equation}\forall i=1,\dots,N, u_{i}:=a_{i}
F_{\mu_{i}}u_{i+1}+b_{i}u_{i+1}F_{\mu_{i}} \textrm{ (except
$b_{N}:=0$)}.\end{equation}

\begin{Lem} The vector $v_1$ in (\ref{hwv}) is a highest weight
vector if and only if \begin{eqnarray}
a_N=\frac{q^{1+\lambda^\mu_N}}{[\lambda^\mu_N]},&\\
(a_{N-1}[\lambda^\mu_{N-1}+1]+b_{N-1}[\lambda^\mu_{N-1}])u_{N}v_0=q^{\lambda^\mu_{N-1}+2}
u_{N}v_0,&\\
(a_{N-1}[\lambda^\mu_{N}]+b_{N-1}[\lambda^\mu_{N}+1])F_{\mu_{N-1}}u_{N+1}v_0=0.&\\
\textrm{ For $i<N-1$: }\qquad\qquad&\\
(a_{i}[\lambda^\mu_{i}+1]+b_{i}[\lambda^\mu_{i}])u_{i+1}v_0=q^{\lambda^\mu_{i}+2}u_{i+1}
v_0,&\\
\left(a_{i}(a_{i+1}[\lambda^\mu_{i+1}+1]+b_{i+1}[\lambda^\mu_{i+1}])\right)F_{\mu_i}u_{i+2}
v_0\ +&\\\left(b_{i}(a_{i+1}[
\lambda^\mu_{i+1}+2]+b_{i+1}[\lambda^\mu_{i+1}+1])\right)F_{\mu_i}u_{i+2}
v_0=0.\label{lambda-1}&\nonumber \end{eqnarray} \end{Lem}

In continuation of the discussion following
Proposition~\ref{6.4}, notice that if $\lambda^\mu_i=0$ then
equation (\ref{lambda-1}) should be stricken, $b_i=0$, and
$a_i=q^2$.

Returning to the general case: If all $\lambda^\mu_i\neq0$:

\begin{eqnarray}
a_N&=&\frac{q^{\lambda^\mu_N+1}}{[\lambda^\mu_N]}.\\\forall
k=1,\dots,N-1:\\ \label{a-n-k}
a_{N-k}&=&q^{\lambda^\mu_{N-k}+2}\frac{[\lambda^\mu_N+\dots+\lambda^\mu_{N-k+1}+k]}{[
\lambda^\mu_N+\dots+\lambda^\mu_{N-k+1}+\lambda^\mu_{N-k}+k]},\\
b_{N-k}&=&-q^{\lambda^\mu_{N-k}+2}\frac{[\lambda^\mu_N+\lambda^\mu_{N-k+1}+k-1]}{[
\lambda^\mu_N+\dots+\lambda^\mu_{N-k+1}+\lambda^\mu_{N-k}+k]}.
\end{eqnarray}

If $\lambda^\mu_{N-1}=0=\dots=\lambda^\mu_{N-R}$ the $a_{N-k}$
just become $q^2$ for $k=1,\dots,R$. This is just the limit of
the equations (\ref{a-n-k}). The corresponding $b_{N-k}=0$
seemingly do not have a nice limit, but recall that instead, we
just cut all branches of the tree marked by $R_{N-i}$,
$i=1.\dots,R$. Actually, in this sense there is a nice limit for
any case in which $\lambda^\mu_i=0$ for some values of
$i=1,\dots, N-1$.

\medskip

\section{One-sidedness} Recall that $\qdet$ is central in
${\mathcal A}^-$.

\begin{Prop}[One-sidedness] If
$\Lambda=(\Lambda_L,\Lambda_R,\lambda)$ and if
\begin{equation}\label{60}
\phi^1_{\Lambda_{n-1},\Lambda_{n}}\circ
\phi^1_{\Lambda_{n-2},\Lambda_{n-1}}\circ\cdots\circ\phi^1_{\Lambda,\Lambda_1}
=\qDet, \end{equation} where $\qDet$ denotes the
operator $M(V_{\Lambda})\ni p\rightarrow \qdet\cdot p \in
M(V_{\Lambda_n})$, then at least one of the pair
$\Lambda_L,\Lambda_R$ is 0. \end{Prop}

We call such a representation {\em one-sided}. We shall see
later that there is a converse to this.

\proof

The proof (sketched) is obtained in 10 installments:

\noindent{\bf 1.} We shall need the following elementary result:
\begin{Lem} Let $a,b\in{\mathbb N}$ with $b\leq a$. Then
$$[a]_q[b]_q=[a+b-1]_q+[a+b-3]_q+\dots +[a-b+1]_q.$$ \end{Lem}

\noindent {\em Proof of Lemma:} Using that
$[a+1]_q=q^{-a}+q^{-a+2}+\dots+q^a$, this follows easily by
counting $q$ exponents. \qed

\smallskip

\noindent{\bf 2.} We have that $\qdet\otimes V\subseteq
{\mathcal A}_1^-\otimes {\mathcal A}_{n-1}^-\otimes V={\mathcal
A}_1^-\otimes ( {\mathcal A}_{n-1}^-\otimes V)=({\mathcal
A}_1^-\otimes {\mathcal A}_{n-1}^-)\otimes V$ and ${\mathcal
A}_{n-1}^-$ is a sum of double tableaux of box size
$(n-1)\times(n-1)$ and similarly ${\mathcal A}_{n}^-$ is a sum
of double tableaux of box size $(n)\times(n)$. By the
Littlewood-Richardson rule, to get $\qdet$ we need to use the
invariant subspace ${\mathcal A}_{n-1}^-(n-1)$ of $(n-1)\times
(n-1)$ minors in ${\mathcal A}_{n-1}^-\otimes V$. We can ignore
contributions from other minors.

\smallskip

\noindent{\bf 3.} We now extend the notation used in
Proposition~\ref{Prop-basis} to also cover the cases of
representations of ${\mathcal U}_q({\mathfrak k}^{\mathbb C}_R)$
in the obvious way. We then have the following extension of said
proposition:

\begin{equation}\label{1}\textrm{If }\quad
v_1=\sum_{k=1,\ell=1}^{i+1,j+1}W_{k,\ell}u_{k,\ell}
\end{equation} is a highest weight vector and $u_{i+1,j+1}=1$,
then

\begin{eqnarray}\label{2}\forall k=1,\dots, i+1,
\ell=1,\dots,j+1:\\u_{k\ell}=\sum_{I_{k,i}\in{\mathcal
J}_{k,i},I_{\ell,j}\in{\mathcal
J}_{\ell,j}}f^\mu(C^\mu(I_{k,i}))f^\nu(C^\nu(I_{\ell,j}))
F^\mu(I_{k,i})F^\nu(I_{\ell,j})v_0.\nonumber \end{eqnarray}

\smallskip

\noindent{\bf 4.} Let ${\mathcal A}_{n-1}^-(n-1)$ denote the
space generated by the $(n-1)\times (n-1)$ minors in ${\mathcal
A}^-$. This is a ${\mathcal U}_q({\mathfrak k}^{\mathbb C})$
module of highest weight
$\Lambda^\mu=(0,0,\dots,0,1)=\Lambda^\nu$. The same kind of
reasoning can be applied to ${\mathcal A}_{n-1}^-(n-1)\otimes
V$. (Notice that ${\mathcal A}_{n-1}^-(n-1)$ is the dual to
${\mathcal A}_{1}^-$.) A highest weight vector $v_0$ in an
irreducible submodule $V_0\subseteq {\mathcal
A}_{n-1}^-(n-1)\otimes \tilde V$ has the form
\begin{equation}v_0= \sum_{k,\ell} A(a+k,b+\ell) \tilde
u_{a+k,b+\ell}\tilde v_0,\label{6} \end{equation} where the
vectors $\tilde u_{a+k,b+\ell}\tilde v_0$, if $k+\ell>0$, have
weights strictly smaller that $\tilde v_0$.

\smallskip

\noindent{\bf 5.} If we insert (\ref{6}) into (\ref{1}) and
isolate the $\tilde v_0$ terms, we get in particular, using
(\ref{60}, (\ref{2}), and since clearly here $(a,b)=(i+1,j+1)$
that \begin{equation} \sum_{k=1}^{i+1}\sum_{\ell=1}^{j+1}
W_{k,\ell}\sum_{I_{k,i}\in{\mathcal
J}_{k,i},I_{\ell,j}\in{\mathcal
J}_{\ell,j}}f^\mu(C^\mu(I_{k,i}))f^\nu(C^\nu(I_{\ell,j}))
F^\mu(I_{k,i})F^\nu(I_{\ell,j})A(i+1,j+1)=\kappa\cdot\qdet
\end{equation} for some constant $\kappa\neq0$. It is easy to
see that $F^\mu(I_{k,i})F^\nu(I_{\ell,j}) A(i+1,j+1))=0$ unless
$(I_{k,i}),I_{\ell,j})=(E_{k,i},E_{\ell,j})$. In the latter case
we get, by (74) in Chapter 1, $(-q^{-1})^{i+j-k-\ell
}A(k,\ell)$.

So \begin{equation} \sum_{k=1}^{i+1}\sum_{\ell=1}^{j+1}
W_{k,\ell}(-q^{-1})^{i+j-k-\ell
}f^\mu(C^\mu(E_{k,i}))f^\nu(C^\nu(E_{\ell,j})A(k,\ell)=\kappa\cdot\qdet.\label{9}
\end{equation}

\smallskip

\noindent{\bf 6.} If both $i+1<n$ and $j-1<n$ we can apply
$F_{\nu_{n-1}}\dots F_{\nu_{j+1}}F_{\mu_{n-1}}\dots
F_{\mu_{i+1}}$ to both sides of (\ref{9}) and get that
$W_{n,n}A(i+1,j+1)=0$; a contradiction.

\smallskip

\noindent{\bf 7.} Let us first assume that $i=j=n$. If we set
$d_{k,\ell}=f^\mu(C^\mu(E_{k,i}))f^\nu(C^\nu(E_{\ell,j})$,
(\ref{9}) becomes

\begin{equation} \sum_{k=1}^{n}\sum_{\ell=1}^{n}
W_{k,\ell}d_{k,\ell}(-q^{-1})^{2n-k-\ell
}A(k,\ell)=\kappa\cdot\qdet.\label{99} \end{equation}

Using (\ref{325}) we can subtract a certain multiple of $\qdet$
in each row such that in the resulting equations

\begin{equation} \sum_{k=0}^{n-1}\sum_{\ell=0}^{n-1}
W_{n-k,n-\ell}b_{k,\ell}A(n-k,n-\ell)=\tilde
\kappa\cdot\qdet,\label{999} \end{equation}we may assume:
$\forall k:b_{k,n}=0$. Of course, this may change the constant
into $\tilde\kappa$. A) If all the remaining $b_{k,\ell}$s are
zero then, naturally, the resulting $\tilde\kappa$ is zero but
that will also imply that each row of the original system
satisfies, up to a constant non-zero multiple, equation
(\ref{325}). In particular, \begin{equation} \sum_{k=1}^{n}
W_{k,n}d_{k,n}(-q^{-1})^{n-k}A(k,n)=\tilde\kappa\cdot\qdet.\label{10}
\end{equation} B) If a non-zero system remains, we can subtract
using column equations (\ref{326}) to remove the terms
$W_{nj}A(n,j)$; $j=1,\dots,n-1$ (the term with $j=n$ has already
been removed. If there still remains an equation
\begin{equation} \sum_{k=0}^{i-1}\sum_{\ell=0}^{j-1}
W_{i-k,j-\ell}b_{k,\ell}A(i-k,j-\ell)=\kappa\cdot\qdet,\label{9a}
\end{equation} we reach a contradiction as in {\bf 6}.

In conclusion:

There is either a column equation

\begin{equation} \sum_{k=1}^{n}
W_{k,n}d_{k,n}(-q^{-1})^{n-k}A(k,n)=\tilde\kappa\cdot\qdet,\label{100}
\end{equation}

or an analogous row  equation 

\begin{equation} \sum_{\ell=1}^{n}
W_{n,\ell}d_{n,\ell}(-q^{-1})^{n-\ell}A(n,\ell)=\tilde\kappa\cdot\qdet.\label{101}
\end{equation}

\medskip

\noindent{\bf 8.} Suppose that we have a row equation
\begin{Lem} \begin{equation} \sum_{\ell=1}^{n}
W_{n,\ell}d_{n,\ell}(-q^{-1})^{n-\ell}A(n,\ell)=\tilde\kappa\cdot\qdet.\label{101}
\end{equation} Then $\lambda^\mu_{n-1}=1$ and $\forall
i=1,\dots,n-2:\lambda^\mu_i=0$. \end{Lem}

\proof We have a PBW basis made up of monomials $W_{n,j_n,
W_{n,j_{n-1}}}, \dots, W_{1,j_1}$. It follows that $\kappa=1$
and it follows from (\ref{101}) and (\ref{325}) that $\forall
k=d_{k,n}=q^{2(n-k)}$. It is easy to see (see {\bf 7.}) that
$d_{k,n}=a_{n-1}a_{n-2}\cdots a_k$. This clearly implies that
$a_{k}=q^2$ for all $k=1,\dots. n-1$.

In particular, $a_{n-1}=q^2$, hence
$$q^2=\frac{q^{1+\lambda^\mu_{n-1}}}{[\lambda^\mu_{n-1}]_q}\Rightarrow
q^{2\lambda^\mu_{n-1}-2}=1\Rightarrow \lambda^\mu_{n-1}=1.$$

Inductively, it follows from (\ref{a-n-k}) that \begin{equation}
q^{\lambda^\mu_{N-k}}\frac{[1+k]}{[\lambda^\mu_{N-k}+1+k]}=1\Rightarrow
\lambda^\mu_{N-k}=0.\end{equation}

\smallskip 

\noindent{\bf 9.} If there is a column equation, it follows in
the same way that $\Lambda_R=(0,0,\dots,0,1)$.

\smallskip

\noindent{\bf 10.} By {\bf 6, 7} what remains are the cases
$i<n, j=n$ and $i=n, j<n$. However, it is clear that they, by
inspection, are covered by the arguments of the case $i=j=n$
simply by eliminating one possibility, so that if $j=n$, we must
have $\Lambda_R=0$ and if $i=n$ we must have $\lambda_L=0$. \qed

\medskip

We have the following converse which is quite straightforward:

\begin{Prop}Let $V_\Lambda=V(\Lambda_L,0,\lambda)$. Set
$\Lambda_0=\Lambda$. Then there exist ${{\mathcal
A}_q^-{\mathcal U}_q({\mathfrak k}^{\mathbb C})}$ intertwining
maps maps $\psi^1_{\Lambda_i,\Lambda_{i+1}}:
V_{\Lambda_{i}}\rightarrow V_{\Lambda_{i+1}} \subset
V_{\Lambda_{i}}\otimes {\mathcal A}^-_q(1)$, for
$i=0,1\dots,n-1$, independent of $\lambda$, such that, with
$\Lambda_n=(\Lambda^0_\mu,0,\lambda+2)$,
$$\psi^1_{\Lambda_{n-1},\Lambda_n}\circ\psi^1_{\Lambda_{n-2},\Lambda_{n-1}}
\circ\cdots\circ\psi^1_{\Lambda,\Lambda_1}=\qDet
.$$ This decomposition is not unique. Furthermore the maps may
be grouped together to form maps of higher degrees, defined by
means of minors of the given degree. \end{Prop}

\bigskip

\section{First order intertwiners}

It is clear that any submodule $V_{\Lambda_1}$ of ${\mathcal
A}_q^-(1)\otimes V_\Lambda$ defines a ${{\mathcal
A}_q^-{\mathcal U}_q({\mathfrak k}^{\mathbb C})}$ equivariant
map $M(V_{\Lambda_1})\rightarrow M(V_\Lambda)$. We shall now see
that there is a unique $\lambda=\lambda(\Lambda_L,\Lambda_R)$
for which this becomes a ${\mathcal U}_q({\mathfrak g}^{\mathbb
C})$ equivariant map. See our forthcoming article for details.
Notice also that the integrality assumption on
$(\Lambda_L,\Lambda_R)$ is not used. \smallskip

We need the following extra information. Modulo ${\mathcal
A}_q^-E_\beta$ it holds:\begin{eqnarray}Z_\beta(W_{i,1})&=&
T_{\mu_{i-1}}T_{\mu_{i-2}}\dots
T_{\mu_{2}}(F_{\mu_1})K_\beta^{-1},\\Z_\beta(W_{i,j})&=&-(q-q^{-1})T_{\nu_{j-1}}T_{\nu_{j-2}}\dots
T_{\nu_{2}}(F_{\nu_1}) T_{\mu_{i-1}}T_{\mu_{i-2}}\dots
T_{\mu_{2}}(F_{\mu_1})K_\beta^{-1} \textrm{ if
}i,j\geq2.\end{eqnarray}

\medskip

\begin{Prop} To any ${\mathcal U}_q({\mathfrak k}^{\mathbb C})$
homomorphism $V_{\Lambda_1}\rightarrow {\mathcal
A}_q^-(1)\otimes V_\Lambda$ there corresponds a unique $\lambda$
such that $\psi_{\Lambda_1,\Lambda}\in Hom_{{\mathcal
U}_q({\mathfrak g}^{\mathbb
C})}(M(V_{\Lambda_1}),M(V_\Lambda))$. \end{Prop}

\medskip

We focus on the case where $\Lambda=(\Lambda_L,0,\lambda)$. Recall (\ref{zbeta}) and  consider \begin{eqnarray}
Z_\beta(W_{N+1}v_0+W_{N}u_Nv_0
+W_{N-1}u_{N-1}v_0+W_{N-2}u_{N-2}v_0 +\dots+W_{1}u_{1}v_0)=\\
q^{-\lambda}\sum_{k=0}^{N-1}T_{\mu_{N-k}}\cdots
T_{\mu_2}T_{\mu_1}(F_{\mu_1})u_{N-k+1}v_0+[\lambda_1+1]u_1v_0=0.\label{sum}
\end{eqnarray}

We may expand the equation into equations for each vector
$F^\mu(I_{1,N})$ in the basis. We claim that the general case
can be reduced by contraction of trees to just the equation for
$F^\mu(E_{1,N})=F_{\mu_1}F_{\mu_2}\cdots F_{\mu_N}v_0$. Here we
get\begin{equation} q^{-\lambda}(1-a_N+a_Na_{N-1}+
a_Na_{N-1}a_{N-2}+\cdots+a_Na_{N-1}a_{N-2}\cdots
a_2+[\lambda+1]a_Na_{N-1}a_{N-2}\cdots a_2 a_1=0. \end{equation}
\begin{eqnarray}
1+a_N&=&q\frac{[\lambda^\mu_N+1]}{[\lambda^\mu_N]},\\
1+a_N+a_Na_{N-1}&=&q^2\frac{[\lambda^\mu_N+1]}{[\lambda^\mu_N]}\frac{[\lambda^\mu_N+\lambda^\mu_
{N-1}+2]}{[\lambda^\mu_N+\lambda^\mu_{N-1}+1]},\\
1+a_N+a_Na_{N-1}+a_Na_{N-1}a_{N-2}&=&q^3
\frac{[\lambda^\mu_N+1]}{[\lambda^\mu_N]}\frac{[
\lambda^\mu_N+\lambda^\mu_{N-1}+2]}{[\lambda^\mu_N+\lambda^\mu_{N-1}+1]}
\frac{[\lambda^\mu_N+\lambda^\mu_{N-1}+\lambda^\mu_{N-2}+
3]}{[\lambda^\mu_N+\lambda^\mu_{N-1}
+\lambda^\mu_{N-2}+2]}.\nonumber\\\end{eqnarray} \begin{eqnarray}
\label{S}S:=1+a_N+a_Na_{N-1}+a_Na_{N-1}a_{N-2}+\cdots+a_Na_{N-1}a_{N-2}\cdots
a_2&=&\\q^{N-1}\frac{[\lambda^\mu_N+1]}{[\lambda^\mu_N]}\cdots
\frac{[\lambda^\mu_N+\cdots+\lambda^\mu_{N-k}+k+1]}{[\lambda^\mu_N+\cdots+\lambda^\mu_{N-k}+k]}
\cdots
\frac{[\lambda^\mu_N+\lambda^\mu_{N-1}+\dots+\lambda^\mu_{2}+
N-1]} {[\lambda^\mu_N+\lambda^\mu_{N-1}+\dots+
\lambda^\mu_{2}+N-2]}. \end{eqnarray}

Comparing to \begin{equation}\label{T}T:=a_Na_{N-1}a_{N-2}\cdots
a_2a_1,\end{equation} one easily obtains
\begin{equation}\label{lambdaeq}
q^{-\lambda}S+[\lambda+1]T=0,\end{equation} which upon divison
becomes \begin{equation}
q^{-\lambda}+[\lambda+1]q^{\lambda^\mu_1+\lambda^\mu_2+\cdots+
\lambda^\mu_N+N}
\frac1{[\lambda^\mu_1+\lambda^\mu_2+\cdots+\lambda^\mu_N+N-1]}=0.
\end{equation} Using the equation $[a+b]=q^{-a}[b]+q^b[a]$, one
easily concludes:
\begin{equation}[\lambda+\lambda^\mu_1+\lambda^\mu_2+\cdots+\lambda^\mu_N+N]=0.
\end{equation} This result can easily be generalized to the
general first order case. It is related to the $q$-Shapovalov
form \cite{mu}.

\end{document}